\documentclass[11pt]{article}
\usepackage{amsmath,amssymb,amsthm,enumerate}
\newtheorem{thm}{Theorem}[section]
\newtheorem{cor}[thm]{Corollary}
\newtheorem{lem}[thm]{Lemma}
\newtheorem{prop}[thm]{Proposition}
\theoremstyle{definition}

\theoremstyle{remark}
\newtheorem*{rem}{Remark}
\newtheorem*{rems}{Remarks}
\numberwithin{equation}{section}

\newcommand{\C}{\mathbb{C}}
\newcommand{\N}{\mathbb{N}}
\newcommand{\R}{\mathbb{R}}
\newcommand{\de}{\delta}
\newcommand{\ep}{\varepsilon}
\newcommand{\ph}{\varphi}
\newcommand{\rd}{\partial}
\newcommand{\hn}{\hat{\nu}}
\newcommand{\calM}{\mathcal{M}}

\newcommand{\ip}[2]{\langle #1 , #2 \rangle}

\DeclareMathOperator{\indx}{index}
\DeclareMathOperator{\nul}{nullity}

\begin{document}

\title{Comparison between Second Variation of Area and \\ 
Second Variation of Energy of a Minimal Surface}

\author{Norio Ejiri \\
Department of Mathematics, 
Faculty of Science and Technology, 
Meijo University \\
1-501 Shiogamaguchi, Tempaku-ku,
Nagoya-shi, Aichi 468-8502 Japan \\
{\sffamily \small ejiri@ccmfs.meijo-u.ac.jp} \and
Mario Micallef \\
Mathematics Institute, University of Warwick, Coventry, CV4 7AL,
U.K. \\ 
{\sffamily \small M.J.Micallef@warwick.ac.uk}} 

\maketitle

\section{Statement and discussion of the results}

The conformal parameterisation of a minimal surface is harmonic 
and therefore, it is natural to compare the Morse index $i_A$ 
of a minimal surface as a critical point of the area functional $A$ 
with its Morse index $i_E$ as a critical point of the energy 
functional $E$.\footnote{Recall that the index of a functional at 
a critical point $F$ is the number, counted with multiplicity, 
of negative eigenvalues of the hessian (i.~e. Jacobi operator) 
of the functional at the critical point. Equivalently, 
the index is the dimension of a maximal subspace of 
the space of infinitesimal variations of $F$ on which 
the hessian is negative definite.} Indeed, 
one way by which minimal surfaces are produced
is by first finding a map which is harmonic with respect to a
fixed conformal structure on the surface and then varying the
conformal structure until we find a harmonic map whose energy is
critical with respect to variations of the conformal structure.
This procedure is well known and has been used very successfully
by Douglas \cite{D}, Courant \cite{C}, Schoen and Yau \cite{ScY},
Sacks and Uhlenbeck \cite{SaU}, Tomi and Tromba \cite{TT} 
and others. 

We now state a theorem relating $i_A$ to $i_E$:

\begin{thm}\label{index_closed} 
Let $F \colon \Sigma_g \to M$ be a (possibly branched) 
minimal immersion of a closed Riemann surface of genus $g$ into 
a Riemannian manifold $M$. Then 
\begin{equation} \label{comp_ineq} 
i_E \leqslant i_A \leqslant i_E + r 
\end{equation} 
where, if $b$ = total number of branch points of $F$ counted with
multiplicity then

\[
r = \begin{cases}6g-6-2b & \text{if }b \leqslant 2g - 3, \\ 
4g - 2 + 2 \left[\frac{-b}{2}\right] & \text{if $2g - 2 \leqslant b 
\leqslant 4g - 4$ and $[x]$ denotes the} \\ 
& \text {largest integer less than or equal to $x$,} \\ 
0 & \text{if $b \geqslant 4g-3$.} \end{cases} 
\]
\end{thm}

Note:  if $g = 0$, $r = 0$ and if $g = 1$, then $r = 2$ if $b = 0$
and $r = 0$ if $b > 0$.\\

\begin{rems} \mbox{} 

\vspace{-3pt} 
\begin{enumerate}[(1)] 
\item $r \leqslant 6g - 6 = \text{real dimension of 
Teichm\"{u}ller space.}$ This is not surprising in light of 
the first paragraph of this paper; see also \S2. 
\item If $g=0$, then any harmonic map is conformal 
(and therefore also minimal) and $i_A = i_E$. 
This is due to the fact that 
the two-sphere carries a unique conformal structure and 
it was proved by the second author in \cite{M2}, Lemma 3.2. 
This paper arose out of that work. 
\item We can also compare the nullity of $F$ 
as a critical point of the area functional with the nullity of $F$ 
as a critical point of the energy functional. 
See Theorem \ref{nul_closed} in \S3 for a precise statement. 
\item Moore has also recently studied the relation between the 
second variations of area and energy in \S5 of \cite{Mo}. 
However, his line of investigation is different from ours. 
In particular, he does not compare the indices of the two functionals. 
\end{enumerate} 
\end{rems} 
Theorem \ref{index_closed} enables us to obtain an upper bound 
on the index of a minimal surface in an arbitrary Riemannian manifold 
which depends on the area and genus of the surface, 
and the dimension and geometry of the ambient manifold; 
see Theorem \ref{gen}. The bound does not depend on 
the second fundamental form of the minimal surface. 

We now consider the index of complete minimal surfaces in 
$\R^d$.\footnote{We shall always mean the area index when 
referring to the index of a minimal surface in $\R^d$.} 
Let $(\Omega_i)_{i \in \N}$ be an exhaustion of a complete minimal 
surface $\Sigma$ by an increasing sequence of compact subsets. 
The index of $\Sigma$ is defined as 
$\sup_{i \to \infty} \indx(\Omega_i)$. 
In her pioneering work \cite{FC}, Fischer-Colbrie showed that a 
complete, oriented minimal surface $\Sigma$ of finite total curvature 
in $\R^3$ has finite index; see also \cite{G}. 
The proof works equally well in $\R^d$, 
$d > 3$ (see, for instance, \cite{N1}). However, 
no bound was given on the index in terms of the total curvature. 
For $d=3$, this was first carried out by Tysk in \cite{T} 
and then improved by Nayatani in \cite{N3}, Theorem 4. 
The case $d > 3$ was treated by the first author in \cite{E1}, 
and also by Cheng and Tysk in \cite{CT}. Among other results, 
they proved that there exists a constant $c_d$, 
depending only on the dimension $d$ of the ambient Euclidean space, 
such that
\[
\indx(\Sigma) \leqslant c_d \int_{\Sigma} (-K)\,dA
\]
where $K$ is the Gauss curvature of $\Sigma$ and 
$dA$ is the element of area on $\Sigma$.\footnote{In \cite{GY} 
Grigory'an and Yau have proved an estimate of this type even for 
a minimal surface in $\R^3$ with boundary.} 
Unfortunately $c_d \to \infty$ as $d \to \infty$. 

The method used in the proof of Theorem \ref{index_closed} 
can also be used to establish the following: 

\begin{thm}\label{index_totcurv} 
Let $\Sigma$ be a complete, oriented, non-planar minimal surface 
with finitely many branch points and of finite total curvature 
in $\R^d$. Then 

\begin{equation}\label{1.1}
\indx(\Sigma) \leqslant \frac{1}{\pi}\int_{\Sigma}(-K)\,dA + 2g - 2 
\leqslant \frac{3}{2\pi}\int_{\Sigma}(-K)\,dA - r + b,
\end{equation}
where $g$ = genus of $\Sigma$, $r$ = number of ends of $\Sigma$ 
and $b$ = total number of branch points counted with multiplicity. 
When $d = 3$ the above inequality may be improved to 

\begin{equation}\label{1.2}
\indx(\Sigma) \leqslant \frac{1}{\pi}\int_{\Sigma}(-K)\,dA + 2g - 3 
\leqslant \frac{3}{2\pi}\int_{\Sigma}(-K)\,dA - r + b -1. 
\end{equation} 

\end{thm}

\begin{rem} 
We can also make statements about $\nul(\Sigma)$; 
see Theorem \ref{nul_totcurv} in \S3 for a precise statement. 
It suffices to state here that we will show that, when $d=3$, 
\[ 
\indx(\Sigma) + \nul(\Sigma) \leqslant 
\frac{1}{\pi}\int_{\Sigma}(-K)\,dA + 2g. 
\] 
This estimate is similar to, but worse than, the one proved 
by Nayatani in \cite{N3}, Theorem 4: 
\begin{equation}\label{1.2nay}
\text{If }\int_{\Sigma}(-K)\,dA \geqslant 8 \pi, \text{ then } 
\indx(\Sigma) + \nul(\Sigma) \leqslant 
\frac{3}{4 \pi}\int_{\Sigma}(-K)\,dA + 3g. 
\end{equation} 
So, perhaps the most striking feature of \eqref{1.1} is that 
it is valid for \emph{all} $d \geqslant 3$. The translational 
Jacobi fields show that, for a non-planar minimal surface in $\R^3$, 
$\nul(\Sigma) \geqslant 3$. In particular, 
\begin{equation}\label{1.2nayindx}
\text{If }\int_{\Sigma}(-K)\,dA \geqslant 8 \pi, \ 
\text{and $d=3$, then } \indx(\Sigma) \leqslant 
\frac{3}{4 \pi}\int_{\Sigma}(-K)\,dA + 3g - 3. 
\end{equation} 
\end{rem}

Fischer-Colbrie showed in \cite{FC} that a complete, oriented 
minimal surface in $\R^3$ of finite index has finite total curvature; 
see also \cite{GL}. 
This does not hold for minimal surfaces in $\R^d$, $d \geqslant 4$ 
because any holomorphic curve (in particular, one with 
infinite total curvature) in $\C^2 = \R^4$ is area-minimizing 
on compact subsets and therefore has index zero. 
(A partial converse to this fact was proved in \cite{M1}.) 
The work of Fischer-Colbrie, of course,
raises the question of obtaining lower bounds on the index of 
a minimal surface in $\R^3$ in terms of the total curvature. 
The first result in this direction was obtained by 
Fischer-Colbrie and Schoen \cite{FCS} (see also \cite{dCP} 
and \cite{Pv}) and states that a complete, stable (i.~e. index zero) 
oriented, minimal surface in $\R^3$ is a plane. Since then, 
several authors have obtained lower bounds, some of which 
we shall compare in the following remarks to 
the upper bounds furnished by \eqref{1.2} and \eqref{1.2nayindx}. 
We refer to \S 7 of \cite{HK} for a fuller discussion 
of the index of minimal surfaces of finite total curvature in $\R^3$. 

\begin{rems}\mbox{} 

\vspace{-3pt} 
\begin{enumerate}[(1)] 
\item If the total curvature of $\Sigma$ in $\R^3$ is $-4\pi$, 
then the Gauss map is 1-1 and therefore $\indx(\Sigma) = 1$. 
It will also follow from Theorem \ref{nul_totcurv} in \S3 that 
the nullity in this case is precisely equal to 3. 
The only complete immersed minimal surfaces with total curvature 
equal to $-4 \pi$ are the catenoid and Enneper's surface; 
see \cite{O}. However, Rosenberg and Toubiana 
have constructed several examples in \cite{RT} 
of branched minimal surfaces in $\R^3$ of total curvature $-4\pi$. 
When $\int_{\Sigma}(-K)\,dA = 4\pi$, genus of $\Sigma = 0$ 
(because the Gauss map is 1-1) and therefore, 
\eqref{1.2} gives $\indx \leqslant 1$. 
This shows that \eqref{1.2} is sharp in this sense. 
We also note that, conversely, L\'{o}pez and Ros proved in \cite{LR} 
that the only complete immersed minimal surfaces in $\R^3$ 
of index 1 are the catenoid and Enneper's surface; see also \cite{MR}. 
\item If $\Sigma_k$ is the Jorge-Meeks $k$-noid of genus zero 
and $k$ ends, then the estimate in \eqref{1.2} can be improved to 
$\indx(\Sigma_k) + \nul(\Sigma_k) \leqslant 2k$; see \S4. 
But in \cite{Ch} it is shown that $\indx(\Sigma_k) \geqslant 2k-3$. 
Therefore, $\indx(\Sigma_k)=2k-3$ and $\nul(\Sigma_k)=3$. 
(Nayatani obtained this result by direct calculation in \cite{N2}; 
see also \cite{MR}.) 
Once again the methods of this paper yield sharp results. 
However, the remarks below indicate that this is not so in general. 
\item If the total curvature of $\Sigma$ in $\R^3$ is $-8 \pi$ 
and $\Sigma$ has genus zero then, \eqref{1.2} and 
the lower bound in \cite{Ch} yield 
$5 \geqslant \indx(\Sigma) \geqslant 3$ whereas according to 
\eqref{1.2nayindx}, $\indx(\Sigma) \leqslant 3$ and therefore, 
$\indx(\Sigma) = 3$. This has been proved by 
the first author and Kotani in \cite{EK}, Corollary 4.3 

The Chen-Gackstatter surface has total curvature equal to $-8 \pi$ 
and genus 1. Montiel and Ros showed that the index of this surface 
is equal to 3 in \cite{MR}, Corollary 15. On the other hand, 
according to \eqref{1.2nayindx} the index is at most 6 and, 
according to \eqref{1.2} the index is at most 7. 
This lack of sharpness of \eqref{1.2} is not surprising as 
it does not take into account any special geometric properties 
of the minimal surface whereas Corollary 15 in \cite{MR} exploits 
the fact that the branching values of the Gauss map of the 
Chen-Gackstatter surface all lie on an equator. 
\item The first author and Kotani \cite{EK} and 
Montiel and Ros \cite{MR} independently proved that 
the index of a complete minimal surface in $\R^3$ of genus zero 
is at most $\frac{1}{2 \pi}\int_{\Sigma}(-K)\,dA - 1$ and that 
generically the index is equal to this number. This shows that 
\eqref{1.2nayindx} (and therefore also \eqref{1.2}) is not sharp 
when the total curvature is less than or equal to $-12 \pi$. 
\item We leave the reader to check that 
\eqref{1.2} and \eqref{1.2nayindx} are also not sharp for 
Bryant's surface and the Hoffman-Meeks surfaces of genus $g$ and 
three ends. 
\end{enumerate}
\end{rems}

This article is organised as follows. In \S \ref{icd} we show that 
the second variation of area for a given normal deformation $s$ 
is less than or equal to the second variation of energy 
for a deformation $v$ whose normal component is $s$. 
Formula \eqref{comp} shows that the difference between 
the two second variations vanishes precisely when 
$v$ is an infinitesimal conformal deformation, 
which is defined by \eqref{2.5}. 
The proofs of the main theorems, based on \eqref{comp} 
and an application of the Riemann-Roch theorem, 
are given in \S \ref{main_proofs}.
In the final section we prove the upper bound on 
the index of a minimal surface in an arbitrary Riemannian manifold 
mentioned above. We also obtain a smaller upper bound 
than that given by Theorem \ref{index_totcurv} for the index of 
a minimal surface of finite total curvature in $\R^3$ 
which has appropriate symmetry; see Theorem \ref{strong_sym}.

\section{Infinitesimal conformal deformations and 
motivation for the proof of Theorem \ref{index_closed}}\label{icd} 

Given a map $F \colon \Sigma \to M$ from a Riemann surface into 
a Riemannian manifold, let $z=x+iy$ be a local complex co-ordinate 
on $\Sigma$. Then,
\begin{equation}\label{2.1}
\text{the energy integrand, $e(F) := 
\frac{1}{2}\left\{|F_x|^2 + |F_y|^2 \right\}$,} 
\end{equation} 
where $F_x = F_*(\rd_x)$, $F_y = F_*(\rd_y)$, 
and the norm $|.|$ is taken with respect to 
the Riemannian metric $\langle\cdot,\cdot\rangle$ on $TM$. The 
\begin{equation}
\text{area integrand, $g(F) := \left\{|F_x|^2|F_y|^2 - \langle
F_x,F_y\rangle^2\right\}^{1/2}$.}
\end{equation}
Therefore
\begin{equation}\label{2.3}
e(F) \geqslant g(F) \text{ with equality if, and only if, 
$F$ is conformal.} 
\end{equation}

We therefore see that, a variation which decreases the energy 
$E := \int_{\Sigma} e(F)\,dxdy$ of a conformal harmonic map $F$ 
must also decrease the area $A := \int_{\Sigma} g(F)\,dxdy$ 
of the map, and therefore $i_A \geq i_E$. Conversely, 
a variation which decreases the area of a conformal harmonic map 
will also decrease the energy if 
we could reparameterise the variation so as to maintain it conformal 
with respect to the initial conformal structure. Of course, 
the obstruction to doing this comes from Teichm\"{u}ller space. 

We now make the above reasoning more formal. 
Let $\nu$ denote the normal bundle of $\Sigma$ and let 
$s \in \Gamma(\nu)$\footnote{$\Gamma$ shall always denote 
the space of smooth sections of a bundle.} be such that
the second variation of area in the direction of $s$, 
$(\de^2 A)(s)$, is negative. Let $\xi$ denote 
the ramified tangent bundle of $\Sigma$, 
i.~e. $\xi$ is the tangent bundle of $\Sigma$ 
twisted at the branch points by an amount equal to 
the order of branching so that $F^*(TM) = \xi \oplus \nu$. 
We wish to find $\sigma_s \in \Gamma(\xi)$ such that

\begin{enumerate}[(1)] 
\item the map $s \mapsto \sigma_s$ is linear,
\item the family of maps corresponding to the variation
vector field $s + \sigma_s$ is a family of conformal maps.
\end{enumerate}
If we succeed, then $(\de^2 E)(s + \sigma_s) = 
(\de^2 A)(s) < 0$. Of course, $\de^2E$ is the hessian 
(or second variation) of the energy functional $E$. 

We will now derive a differential equation for $\sigma_s$ that
will guarantee property $(2)$ at the infinitesimal level. 
Let $I = (-\ep, \ep) \subset \R$, $\ep > 0$ and 
let $V \colon I \times \Sigma \to M$ be such that 
$V(0,\cdot) = F(\cdot)$ and
$(V_*(0,\cdot))(\rd_\ep) = s(\cdot)$. We now want 
$\ph \colon I \times \Sigma \to \Sigma$ such that 
$\ph(0,\cdot)$ = identity and
$\tilde{V}(t,\cdot) := V(t,\ph(t,\cdot))$ is conformal with respect 
to the given fixed conformal structure on $\Sigma$ for all $t \in I$. 
The conformality of $\tilde{V}(t,\cdot)$ can be expressed by
\begin{equation}\label{2.4}
\langle\tilde{V}_*(\rd_z),\tilde{V}_*(\rd_z)\rangle = 0
\end{equation}
where $z$ is a local complex co-ordinate on $\Sigma$ and
$\langle\cdot , \cdot\rangle$ denotes the Riemannian metric on $TM$ 
extended \emph{complex bilinearly} to $T_{\C}M := TM \otimes_{\R} \C$.  
Differentiating \eqref{2.4} with respect to $t$ 
and setting $t=0$ gives: 
\begin{equation}\label{2.5}
\langle\nabla_{\rd_z}(s + \sigma_s),F_z \rangle = 0
\end{equation}
where $\sigma_s = F_*((\phi_*(0,\cdot))\rd_t) \in \Gamma(\xi)$ and 
$\nabla$ is the Levi-Civita connection on $M$ pulled back to $F^*(TM)$ 
and extended complex linearly to $F^*(T_{\C} M)$. 
For obvious reasons a vector field $v \in \Gamma(F^*(TM))$ which 
satisfies $\langle\nabla_{\rd_z}v,F_z \rangle =0$ is called 
an infinitesimal conformal deformation. 

We now recall that the complex structure on $\Sigma$ gives rise to
the splitting $\xi_{\C} := \xi \otimes_{\R} \C = 
\xi^{1,0} \oplus \xi^{0,1}$ where the fibre of $\xi^{0,1}$ 
($\xi^{0,1}$) is locally spanned by $F_z$ ($F_{\bar{z}}$) 
away from the branch points. (The holomorphicity of $F_z$ is required 
to explicitly trivialize $\xi^{1,0}$ on a neighbourhood of 
a branch point.) Therefore we may write 
$\sigma_s = \sigma_s^{1,0} + \sigma_s^{0,1}$. 
Next observe that 
\begin{equation}\label{2.6} 
\langle\nabla_{\rd_z}\sigma_s^{1,0},F_z\rangle = 0 \text{ and } 
\langle\nabla_{\rd_z}\sigma_s^{0,1},F_{\bar{z}}\rangle = 0 
\text{ by conformality of $F$.} 
\end{equation} 
Moreover 
\begin{equation}\label{2.7} 
\ip{\nabla_{\rd_z}s}{F_{\bar{z}}} = 
-\ip{s}{\nabla_{\rd_z}F_{\bar{z}}} = 0 
\text{ by harmonicity of $F$}. 
\end{equation} 
Using \eqref{2.6} and \eqref{2.7} one sees that 
\eqref{2.5} may be re-written as 
\begin{equation}\label{2.8} 
(\nabla_{\rd_z}\sigma_s^{0,1})^{\top} = 
-(\nabla_{\rd_z}s)^{\top} 
\end{equation} 
where the superscript $\top$ denotes orthogonal projection onto 
$\xi_{\C}$. Of course, the global form of \eqref{2.8} is 
\begin{equation}\label{2.9} 
D'\sigma_s^{0,1} = - (\nabla's)^{\top} 
\end{equation} 
where $\nabla' = dz \otimes \nabla_{\rd_z}$, 
$D' =  dz \otimes D_{\rd_{z}}$ and 
$D$ is the connection on $\xi_\C$ induced by $\nabla$ and 
orthogonal projection onto $\xi_{\C}$. 
\eqref{2.9} is the differential equation that 
$\sigma_s$ has to satisfy in order for $v = s + \sigma_s$ 
to be an infinitesimal conformal deformation. 
Theorem \ref{ae_ineq} below essentially asserts the converse.

\begin{thm}\label{ae_ineq} 
Let $F \colon \Sigma \to M$ be a (possibly branched) minimal immersion 
of a closed two-real dimensional oriented surface into 
a Riemannian manifold. For any $s \in \Gamma(\nu)$ and 
$\sigma \in \Gamma(\xi)$ we have 
\begin{equation}\label{3.1} 
(\de^2 A)(s) \leqslant (\de^2 E)(s + \sigma) 
\text{ with equality if and only if $\sigma$ satisfies \eqref{2.9}}. 
\end{equation} 
\end{thm} 

In \eqref{3.1}, $F$ is, of course, being regarded as a harmonic map 
which is conformal (away from the branch points) with respect to 
the conformal structure on $\Sigma$ induced by $F$. 
A more precise relationship between $\de^2 E$ and $\de^2 A$ 
is given by \eqref{comp} below.

\begin{proof} 
One could try to prove this theorem by reversing the argument 
that led to the derivation of \eqref{2.9} but we prefer to give 
a more formal proof that works unchanged also in the case of 
complete minimal surfaces of finite total curvature.

Let $z = x + iy$ be a local complex co-ordinate on $\Sigma$ and 
let $v = s + \sigma$. Then 
\begin{equation}\label{3.2}
  (\de^2 E)(v) = \int_{\Sigma} \left(|\nabla_{\rd_x}v|^2 +
  |\nabla_{\rd_y}v|^2 - \ip{R(v,F_x)F_x}{v} - \ip{R(v,F_y)F_y}{v} 
  \right) \,dxdy 
\end{equation} 
where $R$ is the Riemann curvature tensor of $M$. Now 
\begin{equation}\label{3.3} 
|\nabla_{\rd_x}v|^2 + 
  |\nabla_{\rd_y}v|^2 = 4\,|\nabla_{\rd_z}v|^2. 
\end{equation} 
We let $\perp$ denote orthogonal projection onto 
$\nu_{\C} := \nu \otimes_\R \C$ and obtain
\begin{equation}\label{3.4} 
  \nabla_{\rd_z}v = (\nabla_{\rd_z}v)^{\perp} + \eta + 
  (\nabla_{\rd_z}\sigma^{1,0})^{\top} 
\end{equation} 
where, as suggested by \eqref{2.8}, 
\begin{equation}\label{3.5} 
  \eta := (\nabla_{\rd_z}s)^{\top} + 
  (\nabla_{\rd_z}\sigma^{0,1})^{\top} \, . 
\end{equation} 
On using \eqref{2.6} and \eqref{2.7} in \eqref{3.4} we obtain 
\begin{equation}\label{3.6}
  |\nabla_{\rd_z}v|^2 = |(\nabla_{\rd_z}v)^{\perp}|^2 + |\eta|^2 + 
  |(\nabla_{\rd_z}\sigma^{1,0})^{\top}|^2. 
\end{equation}
Locally, and away from the branch points, we can write
$\sigma^{0,1} = f\,F_{\bar{z}}$ for some locally defined function
$f$. Therefore
\begin{equation}\label{3.7} 
  (\nabla_{\rd_z}\sigma^{0,1})^{\perp} = 0 
  \text{ by harmonicity of $F$}. 
\end{equation} 
\eqref{3.7} allows us to re-write \eqref{3.6} as 
\begin{equation}\label{3.8} 
\begin{split} 
  |\nabla_{\rd_z}v|^2 &= |(\nabla_{\rd_z}s)^{\perp}|^2 + |\eta|^2 + 
  |\nabla_{\rd_z}\sigma^{1,0}|^2 \\ 
  & \qquad + 
  \ip{(\nabla_{\rd_z}s)^{\perp}}{\nabla_{\rd_{\bar{z}}}\sigma^{0,1}} 
  +
  \ip{(\nabla_{\rd_{\bar{z}}}s)^{\perp}}{\nabla_{\rd_z}\sigma^{1,0}}. 
\end{split} 
\end{equation} 
We now calculate the last two terms of \eqref{3.8}: 
\begin{equation}\label{3.9} 
\begin{split} 
\ip{(\nabla_{\rd_z}s)^{\perp}}{\nabla_{\rd_{\bar{z}}}\sigma^{0,1}} 
& = \rd_z \ip{s}{\nabla_{\rd_{\bar{z}}}\sigma^{0,1}} - 
\ip{s}{\nabla_{\rd_z}\nabla_{\rd_{\bar{z}}} \sigma^{0,1}}\, , \\ 
\ip{(\nabla_{\rd_{\bar{z}}}s)^\perp}{\nabla_{\rd_z}\sigma^{1,0}} 
& = \rd_{\bar{z}} \ip{s}{\nabla_{\rd_z}\sigma^{1,0}} - 
\ip{s}{\nabla_{\rd_{\bar{z}}}\nabla_{\rd_z} \sigma^{1,0}} \, . 
\end{split} 
\end{equation} 
But, from \eqref{3.7} and \eqref{3.5} we have 
$\nabla_{\rd_z}\sigma^{0,1} = (\nabla_{\rd_z}\sigma^{0,1})^{\top} = 
\eta - (\nabla_{\rd_z}s)^{\top}$ and therefore 
\begin{equation}\label{3.10} 
\begin{split} 
  \ip{\nabla_{\rd_z}\nabla_{\rd_{\bar{z}}}\sigma^{0,1}}{s} 
  &= \ip{R(F_z,F_{\bar{z}})\sigma^{0,1}}{s} + 
  \ip{\nabla_{\rd_{\bar{z}}}(\eta - (\nabla_{\rd_z}s)^{\top})}{s} \\ 
  &= \ip{R(F_z,F_{\bar{z}})\sigma^{0,1}}{s} + 
  |(\nabla_{\rd_z}s)^{\top}|^2 - 
  \ip{\eta}{(\nabla_{\rd_{\bar{z}}}s)^{\top}}. 
\end{split} 
\end{equation} 
Similarly, 
\begin{equation}\label{3.11} 
\ip{\nabla_{\rd_{\bar{z}}}\nabla_{\rd_z}\sigma^{1,0}}{s} = 
\ip{R(F_{\bar{z}},F_z)\sigma^{1,0}}{s} + 
  |(\nabla_{\rd_{\bar{z}}}s)^\top|^2 - 
  \ip{\bar{\eta}}{(\nabla_{\rd_z}s)^{\top}}. 
\end{equation} 
Taking \eqref{3.9}, \eqref{3.10} and \eqref{3.11} into account in 
\eqref{3.8}, integrating and using Stokes's theorem gives 
\begin{equation}\label{3.12} 
\begin{split} 
  \int_{\Sigma} |\nabla_{\rd_z}v|^2 \,dxdy & = \int_{\Sigma} 
  \left(|(\nabla_{\rd_z}s)^{\perp}|^2 + |\eta|^2 + 
  |\nabla_{\rd_z}\sigma^{1,0}|^2 - 
  2\,|(\nabla_{\rd_z}s)^{\top}|^2 \right. \\ 
  & \qquad \quad - \ip{R(F_z,F_{\bar{z}})\sigma^{0,1}}{s} - 
  \ip{R(F_{\bar{z}},F_z)\sigma^{1,0}}{s} \\ 
  & \left. \qquad \quad + \ip{\eta}{(\nabla_{\rd_{\bar{z}}}s)^{\top}} 
  + \ip{\bar{\eta}}{(\nabla_{\rd_z}s)^{\top}}\right ) \,dxdy. 
\end{split} 
\end{equation} 
We now deal with the last two terms in \eqref{3.2}: 
\begin{equation}\label{3.13} 
  \begin{split} 
  \ip{R(v,F_x)F_x}{v} + \ip{R(v,F_y)F_y}{v} 
  & = 4\,\ip{R(v,F_z)F_{\bar{z}}}{v} \\ 
  & = 4\, \left (\ip{R(s,F_z)F_{\bar{z}}}{s} + 
  \ip{R(\sigma^{0,1},F_z)F_{\bar{z}}}{\sigma^{1,0}} \right. \\ 
  & \left. \qquad + \ip{R(s,F_z)F_{\bar{z}}}{\sigma^{1,0}} + 
  \ip{R(\sigma^{0,1},F_z)F_{\bar{z}}}{s} \right ). 
  \end{split} 
\end{equation} 
By the second Bianchi identity, 
\begin{equation}\label{3.14} 
  \begin{split} 
\ip{R(\sigma^{0,1},F_z)F_{\bar{z}}}{s} + 
\ip{R(F_z,F_{\bar{z}})\sigma^{0,1}}{s} &= 0, \\ 
\ip{R(\sigma^{1,0},F_{\bar{z}})F_z}{s} + 
\ip{R(F_{\bar{z}},F_z)\sigma^{1,0}}{s} &= 0. 
  \end{split} 
\end{equation} 
Using \eqref{3.3}, \eqref{3.12}, \eqref{3.13} and \eqref{3.14} in 
\eqref{3.2} yields: 
\begin{equation}\label{3.15} 
  \begin{split} 
  (\de^2 E)(v) &= 4 \int_{\Sigma} \left( 
  |(\nabla_{\rd_z}s)^{\perp}|^2 + |\nabla_{\rd_z}\sigma^{1,0}|^2 + 
  |\eta|^2 - 2\,|(\nabla_{\rd_z}s)^{\top}|^2 \right. \\ 
  & \qquad \qquad - \ip{R(s,F_z)F_{\bar{z}}}{s} - 
  \ip{R(\sigma^{0,1},F_z)F_{\bar{z}}}{\sigma^{1,0}} \\ 
  & \left. \qquad \qquad + \ip{\eta}{(\nabla_{\rd_{\bar{z}}}s)^{\top}} 
  + \ip{\bar{\eta}}{(\nabla_{\rd_z}s)^{\top}} \right) \,dxdy. 
  \end{split} 
\end{equation} 
An integration by parts shows that 
\begin{equation}\label{3.16} 
  \int_{\Sigma} |\nabla_{\rd_z}\sigma^{1,0}|^2\,dxdy = \int_{\Sigma} 
  |\nabla_{\rd_{\bar{z}}}\sigma^{1,0}|^2 \,dxdy + \int_{\Sigma} 
  \ip{R(F_z,F_{\bar{z}})\sigma^{1,0}}{\sigma^{0,1}}\,dxdy 
\end{equation} 
which, together with \eqref{3.5} and 
the second Bianchi identity, gives: 
\begin{equation}\label{3.17} 
  \begin{split} 
  \int_{\Sigma} |\nabla_{\rd_z}\sigma^{1,0}|^2\,dxdy 
  &= \int_{\Sigma} \left( |\eta|^2 + |(\nabla_{\rd_z}s)^{\top}|^2 - 
  \ip{\eta}{(\nabla_{\rd_{\bar{z}}}s)^{\top}} \right. \\ 
  & \left. \qquad - \ip{\bar{\eta}}{(\nabla_{\rd_z}s)^{\top}} + 
  \ip{R(\sigma^{1,0},F_{\bar{z}})F_z}{\sigma^{0,1}}\right) \,dxdy. 
  \end{split} 
\end{equation} 
On substituting \eqref{3.17} in \eqref{3.15} we obtain 
\begin{equation}\label{comp} 
\begin{split} 
(\de^2 E)(v) &=4 \int_{\Sigma} \left(
|(\nabla_{\rd_z}s)^{\perp}|^2 - |(\nabla_{\rd_z}s)^{\top}|^2 - 
\ip{R(s,F_z)F_{\bar{z}}}{s} + 2|\eta|^2 \right)\,dxdy \\ 
&= (\de^2 A)(s) + 8 \int_{\Sigma} |\eta|^2\,dxdy \, . 
\end{split} 
\end{equation} 
The proof of the theorem is complete. 
\end{proof} 

\section{Proof of Theorems \ref{index_closed} and 
\ref{index_totcurv}}\label{main_proofs} 

\begin{proof}[Proof of Theorem \ref{index_closed}]
The inequality $i_A \geqslant i_E$ follows immediately from 
Theorem \ref{ae_ineq}. For the other inequality, 
let $S$ be a maximal subspace on which $\de^2 A < 0$. 
By the Fredholm alternative, we may solve \eqref{2.9} if, 
and only if, $(\nabla_{\rd_z}s)^\top$ is 
orthogonal to $\ker D'{^*}$, where $D'{^*}$ is the adjoint of $D'$. 
Now an integration by parts shows that $D'{^*} = i*\bar{\rd} \colon 
\Gamma(\xi^{0,1} \otimes \kappa) \to \Gamma(\xi^{0,1})$ where 
$\kappa$ is the line bundle of holomorphic one-forms on $\Sigma$ 
and $*$ is the Hodge star operator. 
Therefore $\ker D'{^*} = H^0(\xi^{0,1} \otimes \kappa)$ = 
space of holomorphic sections of $\xi^{0,1} \otimes \kappa$. Let
$h^0(\xi^{0,1} \otimes \kappa)$ = complex dimension of
$H^0(\xi^{0,1} \otimes \kappa)$. Then, we may find a subspace
$\tilde{S} \subset S$ of real dimension $\geqslant$ 
$\dim S - 2h^0(\xi^{0,1} \otimes \kappa)$ for which \eqref{2.9} 
has a solution $\sigma_s$ whenever $s \in \tilde{S}$. 
Moreover, we may arrange for $\sigma_s$ to depend linearly on $s$, 
since the equation for $\sigma_s$ is linear. 
Let $\hat{S} = \{ s + \sigma_s \mid s \in \tilde{S} \} 
\subset \Gamma(F^*(TM))$. Then, by Lemma \ref{ae_ineq}, 
$\de^2 E \vert_{\hat{S}} < 0$ and therefore, 
$i_E \geqslant \dim \hat{S} = \dim \tilde{S} \geqslant i_A - r$ 
where $r = 2 h^0(\xi^{0,1} \otimes \kappa)$. We now calculate 
$h^0(\xi^{0,1} \otimes \kappa)$: $c_1(\xi^{0,1}) = 2g-2-b$ and 
$c_1(\kappa) = 2g-2$ and therefore, by the theorem of 
Riemann-Roch, we have 
\[ 
h^0(\xi^{0,1}\otimes \kappa) = 3g-3-b + h^0(\xi^{1,0}). 
\] 
If $b \leqslant 2g-3$, $c_1(\xi^{1,0}) < 0$ and 
$h^0(\xi^{0,1}\otimes \kappa) = 3g - 3 - b$. 
If $b \geqslant 4g - 3$, $c_1(\xi^{0,1}\otimes \kappa) < 0$ and  
$h^0(\xi^{0,1} \otimes \kappa) = 0$. 
If $2g - 2 \leqslant b \leqslant 4g - 4$, then 
$0 \leqslant c_1(\xi^{0,1} \otimes \kappa) \leqslant 2g - 2$ and 
therefore, by Clifford's theorem (see, for example, \cite{GH}), 
$h^0(\xi^{0,1} \otimes \kappa) \leqslant 
\left[\frac{4g - 2 - b}{2}\right]$. 

The proof of Theorem \ref{index_closed} is complete. 
\end{proof} 

Recall that the nullity $n$ of a functional at a critical point $F$ is 
the dimension of the space of Jacobi fields of the functional at $F$. 
If the index of $F$ is $i$ then $i+n$ is equal to 
the dimension of a maximal subspace of the space of 
infinitesimal variations of $F$ on which the hessian of the functional 
is negative semidefinite. A minor modification 
(which will be left to the reader) of the proof of 
Theorem \ref{index_closed} yields: 

\begin{thm}\label{nul_closed} 
Let $F \colon \Sigma_g \to M$ and $r$ be as in 
Theorem \ref{index_closed} and let 
\begin{align*} 
n_A &= \text{nullity of $F$ as 
a critical point of the area functional $A$,} \\ 
n_E&= \text{nullity of $F$ as 
a critical point of the energy functional $E$,} \\ 
{n_E}^T&= \text{dimension of the space of \emph{purely tangential}} \\ 
&\hphantom{\mbox{}=\mbox{}} \text{Jacobi fields of $F$, 
as a critical point of $E$.} 
\end{align*} 
Then 
\[  
i_E + n_E - n_E^T \leqslant i_A + n_A \leqslant i_E + n_E - n_E^T + r. 
\] 
\end{thm} 
The following comparison of the nullities of energy and area 
follows immediately from the inequalities in 
Theorem \ref{index_closed} and Theorem \ref{nul_closed}. 
\[ 
n_E - n_E^T - r \leqslant n_A \leqslant n_E - n_E^T + r. 
\] 

We now move on to the proof of Theorem \ref{index_totcurv}. 
First, we recall that (see, for example, \cite{L}) 
if a minimal surface $\Sigma$ in $\R^d$ 
has finite total curvature and finitely many branch points, 
then $\Sigma$ is conformally diffeomorphic to a closed Riemann surface 
$\bar{\Sigma}$ with finitely many punctures $\{p_1, \dotsc ,p_k\}$ 
corresponding to the ends of $\Sigma$. Recall, too, that $G_{2,d}$, 
the Grassmannian of oriented two-planes in $\R^d$, 
may be identified with the quadric $Q_{d-2} \subset \C P^{d-1}$ 
defined by $\{[z] \mid z = (z_1, \dotsc , z_d) 
\in \C^d \setminus \{0\}, \ \sum_{i=1}^d (z_i)^2 = 0\}.$ 
Furthermore, the Gauss map $G \colon \Sigma \to G_{2,d} = Q_{d-2}$ 
is holomorphic and extends to a holomorphic map 
$\bar{G} \colon \bar{\Sigma} \to G_{2,d}$. 
Let $\gamma$ be the tautological two-plane bundle over $G_{2,d}$ 
and let $\bar{\xi} = \bar{G}^*(\gamma)$. Then 
$\left. \bar{\xi} \right|_{\Sigma} = \xi$, 
the ramified tangent bundle of $\Sigma$ and 
$c_1(\bar{\xi}) = \frac{1}{2 \pi} \int_{\Sigma} K \, dA$. 
Similarly, let $\bar{\nu}$ be the orthogonal complement of $\bar{\xi}$ 
in $\bar{\Sigma} \times \R^d$ and then 
$\left. \bar{\nu} \right|_{\Sigma} = \nu$. 
The hessians $\de^2 E$ and $\de^2 A$ both extend to 
sections of $\bar{\Sigma} \times \R^d$ and $\bar{\nu}$ respectively 
and, in \cite{FC}, \cite{G} and \cite{N1} it is shown that 
$i_A(\Sigma) = i_A(\bar{\Sigma}).$ We also have 
$i_E(\bar{\Sigma}) = 0$ and $n_E(\bar{\Sigma}) = d$ 
because the tangent bundle of $\R^d$ is trivial and 
the Levi-Civita connection on $\R^d$ is simply the exterior derivative. 
Finally, the projection of a constant vector field in $\R^d$ 
onto $\bar{\nu}$ is a Jacobi field of the area functional and therefore, 
if $\Sigma$ is non-planar then $n_A(\Sigma) \geqslant d$. We shall, 
in fact, prove the following more precise theorem. 

\begin{thm}\label{nul_totcurv} 
Let $\Sigma$ be a complete, oriented, non-planar minimal surface 
in $\R^d$ with finitely many branch points and 
of finite total curvature. Define $g$, $b$ and $r$ as in 
Theorem \ref{index_totcurv}. Then 
\begin{equation}\label{totcurv} 
\indx(\Sigma) + \nul(\Sigma) \leqslant 
\frac{1}{\pi} \int_{\Sigma} (-K) \, dA + 2g - 2 + d \leqslant 
\frac{3}{2 \pi} \int_{\Sigma} (-K) \, dA - r + b + d. 
\end{equation} 
Furthermore, if $d=3$ then \eqref{totcurv} can be improved to 
\begin{equation}\label{totcurv_3} 
\indx(\Sigma) + \nul(\Sigma) \leqslant 
\frac{1}{\pi} \int_{\Sigma} (-K) \, dA + 2g \leqslant 
\frac{3}{2 \pi} \int_{\Sigma} (-K) \, dA - r + b + 2. 
\end{equation} 
\end{thm} 
\begin{rem} 
Inequalites \eqref{1.1} and \eqref{1.2} in Theorem \ref{index_totcurv} 
result from using, in \eqref{totcurv} and \eqref{totcurv_3}, 
the observation that $\nul(\Sigma) \geqslant d$ 
for a non-planar minimal surface. 
\end{rem} 

\begin{proof} 
It is clear that the proof of Theorem \ref{ae_ineq} also works for 
$\de^2 E$ and $\de^2 A$ extended to $\bar{\Sigma}$. Furthermore, 
\[ 
c_1(\bar{\xi}^{1,0}) = c_1(\bar{\xi}) = 
\frac{1}{2 \pi}\int_{\Sigma} K \, dA < 0. 
\] 
Therefore, by the Riemann-Roch theorem we have 
\[ 
h^0(\bar{\xi}^{0,1} \otimes \kappa) = 
- c_1(\bar{\xi}) + c_1(\kappa) + 1 - g = 
\frac{1}{2 \pi}\int_{\Sigma} (-K) \, dA + g - 1. 
\] 
The proof of the first inequality in \eqref{totcurv} 
can now be completed by an argument identical to that of 
the proof of Theorem \ref{nul_closed}. The inequality 
\[ 
\frac{1}{2 \pi}\int_{\Sigma} (-K) \, dA \geqslant 
2g - 2 + r - b 
\] 
is due to Osserman; see for instance, \cite{L}. 

The proof of \eqref{totcurv_3} when $d=3$ requires the following lemma. 
\begin{lem} \label{hol} 
Let $F \colon \Sigma \to \R^3$ be a generalised minimal immersion 
of a Riemann surface $\Sigma$. Let $\hn$ be 
a smooth unit normal vector field on $\Sigma$. 
(Such a section of $\nu$ exists because $\Sigma$ is orientable.) 
Then $(\rd \hn)^{\top} \in H^0(\xi^{0,1} \otimes \kappa)$. 
\end{lem} 
\begin{proof} 
By \eqref{2.7} we have that 
\[ 
(\rd \hn)^{\top} = \ip{\rd_z \hn}{F_z} 
\frac{F_{\bar{z}}}{|F_z|^2} \otimes dz = 
- \ip{\hn}{F_{zz}} \frac{F_{\bar{z}}}{|F_z|^2} \otimes dz . 
\] 
The conformality and harmonicity of $F$ then show that 
$D''\big((\rd \hn)^{\top}\big) = 0$, where 
$D'' := d \bar{z} \otimes D_{\rd_{\bar{z}}}$, 
i.~e., $(\rd \hn)^{\top}$ is a holomorphic section of 
$\xi^{0,1} \otimes \kappa$, as claimed. 
\end{proof} 
We now return to the proof of \eqref{totcurv_3}. Denote by $\calM$ 
the space of meromorphic functions on $\bar{\Sigma}$. 
By Lemma \ref{hol}, 
\[ 
H^0(\xi^{0,1} \otimes \kappa) = \{g(\rd \hn)^{\top} \mid 
g \in \calM, \ [g] + [(\rd \hn)^{\top}] \geqslant 0\}, 
\] 
where $[\cdot]$ denotes divisor. It is convenient to define 
\begin{equation} \label{ml} 
\calM_L := \{g \in \calM \mid [g] + [(\rd \hn)^{\top}] \geqslant 0\}, 
\end{equation} 
where, of course, $L := \xi^{0,1} \otimes \kappa$. 
If $s \in \Gamma(\nu)$, then $s = f \hn$ for some 
$f \in C^{\infty}(\bar{\Sigma})$. 
Therefore, the Fredholm alternative\footnote{as in 
the proof of Theorem \ref{index_closed}} 
for the solvability of \eqref{2.9} can now be stated as 
the following condition on $f$: 
\begin{equation}\label{fred} 
\int_{\Sigma}f \bar{g} K \, dA = 0 \ \forall \, g \in \calM_L, 
\end{equation} 
where we have used $|(\rd \hn)^{\top}|^2 = -K$. 
Since the constant function 1 belongs to $\calM_L$, we see that 
the $\R$-codimension of the space of real valued smooth functions $f$ 
satisfying \eqref{fred} is $2 h^0(\xi^{0,1} \otimes \kappa) - 1 = 
\frac{1}{\pi} \int_{\Sigma} (-K) \, dA + 2g - 3.$ 
\end{proof} 

\section{More area index estimates} \label{particular} 
Theorem \ref{index_closed} provides an upper bound on $i_A$ whenever 
$i_E$ may be estimated, and it is often easier 
to estimate $i_E$ than $i_A$ because 
$\de^2 E$ does not involve the second fundamental form. 
For instance, if $M$ has nonpositive sectional curvature 
then $i_E$ is easily seen to be zero from \eqref{3.2}, 
thereby yielding: 
\begin{cor} \label{nonpos}
Let $\Sigma$ be a closed minimal surface in a Riemannian manifold $M$ 
of nonpositive sectional curvature. Then $i_A \leqslant r 
\leqslant 6g - 6$ where $r$ is as in Theorem \ref{index_closed}. 
\end{cor} 
In particular, the index of a closed minimal surface of genus $g$ 
in a flat torus is at most $6g-6$. Corollary \ref{nonpos} 
has already been noted in \cite{E3} by considering 
an energy functional on Teichm\"{u}ller space. 

Theorem \ref{index_closed} can be refined when the ambient space is 
3-dimensional in a manner similar to that
in Theorem \ref{index_totcurv} for minimal surfaces 
of finite total curvature in $\R^3$. More precisely, 

\begin{thm}\label{index_closed3} 
Let $F \colon \Sigma_g \to M^3$ be a (possibly branched) 
minimal immersion of a closed Riemann surface of genus $g$ into 
a three-dimensional space-form $M^3$. 
If $F$ is not totally geodesic then 
\begin{equation} \label{comp_ineq3} 
i_E \leqslant i_A \leqslant i_E + r - 1
\end{equation} 
where $r$ is as in Theorem \ref{index_closed}. 
\end{thm} 

The proof consists in observing that Lemma \ref{hol} 
is still valid in this setting. The argument then proceeds 
as in the proof of \eqref{totcurv_3}. 

The Clifford torus in $\R P^3$ is stable as a harmonic map. 
(In \cite{E2}, the first author has determined all harmonic tori 
in $\R P^3$ that minimize energy in their homotopy class.) 
However, the Clifford torus is unstable as a minimal surface 
and it is not hard to show that its index is 1. Furthermore 
$r = 2$. This is a situation in which $i_A < i_E$ and yet 
\eqref{comp_ineq3} is still sharp. 

For a general minimal immersion we can prove 

\begin{thm}\label{gen} 
Let $F \colon \Sigma_g \to M^d$ be a (possibly branched) 
minimal immersion of a closed Riemann surface of genus $g$ into 
a Riemannian manifold $M$. Then 
\begin{equation}\label{eq_gen_h} 
i_E + n_E \leqslant d \, C(M) \, \text{Area}\,(F(\Sigma_g)) 
\end{equation} 
and therefore, 
\begin{equation}\label{eq_gen_m} 
i_A + n_A \leqslant d \, C(M) \, \text{Area}\,(F(\Sigma_g)) + r 
\end{equation} 
where $r$ is as in Theorem \ref{index_closed} and $C(M)$ 
is a constant which depends on the second fundamental form of 
an isometric embedding of $M$ into Euclidean space. 
\end{thm} 
\begin{rems} 
The main interest in \eqref{eq_gen_m} is that 
no assumption is made on the second fundamental form of $F$. 
Colding and Minicozzi showed (Theorem 1.1 in \cite{CM}) that, 
given $A>0$ and a positive integer $g$, 
there are at most finitely many closed embedded minimal surfaces 
of genus $g$ and with area at most $A$ 
in a closed orientable 3-manifold with a bumpy metric. 
In particular, the Morse index of an embedded minimal surface 
in a 3-manifold with a bumpy metric is bounded 
by its genus and its area but no explicit bound like 
\eqref{eq_gen_m} is given in \cite{CM}. 
\end{rems} 
\begin{proof} 
Let $h(t)$ be the trace of the heat kernel of $\Sigma_g$ with 
the metric induced by $F$. Then, since $F$ is an isometric harmonic map, 
Proposition 2.2 in \cite{U} asserts that 
\begin{equation}\label{eq_gen_hu} 
i_E + n_E \leqslant d \, \inf_{t>0}e^{2 \mu t}h(t)
\end{equation} 
where $\mu$ is an upper bound of the sectional curvatures of $M^d$. 

Now consider $M$ to be isometrically embedded in some Euclidean space. 
The method of proof of Theorem 5 on pages 991-993 in \cite{CT} 
can now be employed to obtain 
\begin{equation}\label{h_bound} 
h(t) \leqslant \frac{\alpha_1}{(1 - e^{-\alpha_2}t)^2} \, 
\text{Area}\,(F(\Sigma_g)) 
\end{equation} 
where $\alpha_1$ and $\alpha_2$ are constants which depend 
on the second fundamental form of the isometric embedding of $M$ 
into Euclidean space. The bound \eqref{eq_gen_h} now follows 
from \eqref{eq_gen_hu} and \eqref{h_bound}. 
\end{proof} 
The indexes of some minimal surfaces of finite total curvature 
have already been mentioned in the remarks at the end of \S1. 
Many known examples of minimal surfaces in $\R^3$ which are embedded 
or, at least, whose ends are embedded, have symmetries 
which allow refinements to Theorem \ref{index_totcurv}. 
The next proposition makes precise the type of symmetry that 
the minimal surface is required to have. It includes the notion 
of strong symmetry with respect to a plane introduced by 
Cos\'{i}n and Ros in \cite{CR}, Definition 1; see also Lemma 4 
of the same article. 

\begin{prop} \label{sym} 
Let $\Sigma$ be a Riemann surface and let $F \colon \Sigma \to \R^3$ 
be a generalised minimal immersion which is complete and 
of finite total curvature. Suppose there exists an isometry 
$\Theta \colon \R^3 \to \R^3$ and a diffeomorphism 
$\theta \colon \Sigma \to \Sigma$ such that 
$\Theta \circ F = F \circ \theta$. 
Let $\calM_L$ be defined by \eqref{ml}. 
If $\theta$ is anti-holomorphic and $g \in \calM_L$ 
then $\overline{g \circ \theta} \in \calM_L$. 
If $\theta$ is holomorphic and $g \in \calM_L$ 
then $g \circ \theta \in \calM_L$. 
\end{prop} 

\begin{proof} 
We shall only give the proof when $\theta$ is anti-holomorphic; 
the proof when $\theta$ is holomorphic is similar and, in fact, 
much easier. 

From $F(\theta(z)) = \Theta(F(z))$ and $\rd \theta = 0$ 
we obtain: 
\[ 
\frac{\rd F}{\rd \bar{w}} (\theta(z)) 
\frac{\rd \bar{\theta}}{\rd z} = 
\Theta_0 \left( \frac{\rd F}{\rd z}(z)\right) 
\] 
where $\Theta_0 \in SO(3)$ is the non-translational part of $\Theta$; 
of course, $w$ is a local complex coordinate defined on 
a neighbourhood of $\theta(z)$. It follows that 
\begin{equation} \label{hnu_sym}
\hn(\theta(z)) = \begin{cases} 
\Theta_0(\hn(z)), & \text{if $\det \Theta_0 = -1$,} \\ 
-\Theta_0(\hn(z)), & \text{if $\det \Theta_0 = 1$.}
\end{cases} 
\end{equation}
Differentiating \eqref{hnu_sym} yields: 
\[ 
\frac{\rd \hn}{\rd \bar{w}} (\theta(z)) 
\frac{\rd \bar{\theta}}{\rd z} (z) = 
\pm \Theta_0 \left( \frac{\rd \hn}{\rd z}(z)\right), 
\text{ i.e., } \theta^*\big( (\bar{\rd} \hn)^{\top} \big) = 
\pm \Theta_0 \big( (\rd \hn)^{\top} \big) . 
\] 
Therefore, $(\rd \hn)^{\top}$ has a zero of order $Q$ at $z$ if, 
and only if, it also has a zero of order $Q$ at $\theta(z)$. 
The proposition follows immediately. 
\end{proof} 

\begin{thm} \label{strong_sym}
Let $\Sigma$ be a Riemann surface and let $F \colon \Sigma \to \R^3$ 
be a generalised minimal immersion which is complete and 
of finite total curvature. Suppose there exists an isometry 
$\Theta \colon \R^3 \to \R^3$ and an anti-holomorphic involution  
$\theta \colon \Sigma \to \Sigma$ such that 
$\Theta \circ F = F \circ \theta$. Then 
\begin{equation}\label{totcurv_3sym} 
\indx(\Sigma) + \nul(\Sigma) \leqslant 
\frac{1}{2 \pi} \int_{\Sigma} (-K) \, dA + g + 2. 
\end{equation} 
In particular, since $\nul(\Sigma) \geqslant 3$, we have 
\begin{equation}\label{indx_3sym} 
\indx(\Sigma) \leqslant 
\frac{1}{2 \pi} \int_{\Sigma} (-K) \, dA + g - 1. 
\end{equation} 
\end{thm} 

\begin{proof} 
Proposition \ref{sym} enables us to define 
$\rho \colon \calM_L \to \calM_L$ by 
$\rho (g) := \overline{g \circ \theta}$, where 
$\calM_L$ is defined by \eqref{ml}. Then $\rho$ is linear and 
$\rho^2 = \text{identity}$. Therefore, $\calM_L = \calM_L^+ \oplus 
\calM_L^-$ where $\calM_L^+$ and $\calM_L^-$ are respectively the 
$+1$ and $-1$ eigenspaces of $\rho$. Similarly, we can write 
$C^{\infty}(\bar{\Sigma}) = C^{\infty}_+(\bar{\Sigma}) \oplus 
C^{\infty}_-(\bar{\Sigma})$. 

Next observe that if $f \in C^{\infty}(\bar{\Sigma})$ then 
$(\de^2 A)(f \hn) = (\de^2 A)(f_+ \hn) + (\de^2 A)(f_- \hn)$, where 
$f_{\pm} := \frac12 (f \pm f \circ \theta)$. Now let $S_+$ be 
a maximal subspace of $C^{\infty}_+(\bar{\Sigma})$ on which 
$\de^2 A \leqslant 0$ and define $S_-$ similarly. 
A moment's thought will reveal that $S := S_+ \oplus S_-$ is then 
a maximal subspace of $C^{\infty}(\bar{\Sigma})$ on which 
$\de^2 A \leqslant 0$. 

Let $\{f_1, \dotsc , f_p\}$ and $\{f_{p+1}, \dotsc , f_q \}$ 
be bases of $S_+$ and $S_-$ respectively and let 
$\{g_1, \dotsc , g_{\mu}\}$ and $\{g_{\mu + 1}, \dotsc , g_{\nu}\}$ 
be bases of $\calM_L^+$ and $\calM_L^-$ respectively. Then, 
for $j \in \{1, \dotsc , q \}$ and $\alpha \in \{1, \dotsc , \nu \}$ 
we have: 
\begin{align*} 
\int_{\Sigma}f_j \bar{g}_{\alpha} K \, dA &= 
\int_{\Sigma}(f_j \circ \theta) (\bar{g}_{\alpha} \circ \theta ) 
\theta^*(K \, dA) \\ 
&= \pm \int_{\Sigma}f_j \bar{g}_{\alpha} K \, dA . 
\end{align*} 
Therefore, $\int_{\Sigma}f_j \bar{g}_{\alpha} K \, dA $ is either 
real or purely imaginary. It follows that the $\R$-codimension of 
the space of real valued smooth functions $f$ 
satisfying \eqref{fred} is $ h^0(\xi^{0,1} \otimes \kappa) = 
\frac{1}{2 \pi} \int_{\Sigma} (-K) \, dA + g - 1.$ 
\eqref{totcurv_3sym} and \eqref{indx_3sym} follow immediately. 
\end{proof} 

The Jorge-Meeks $k$-noid $\Sigma_k$ of genus zero and $k$ ends 
has total curvature equal to $4 \pi (k - 1)$ and 
is strongly symmetric in the sense of Cos\'{i}n and Ros in \cite{CR}. 
Therefore, we may apply Theorem \ref{strong_sym} to conclude, 
as asserted in a remark in \S1, that 
$\indx(\Sigma_k) + \nul(\Sigma_k) \leqslant 2k.$ 
 
\end{document}